\documentclass{article}
\usepackage{amsfonts}
\usepackage{marvosym}
\title{A Proof of The Lonely Runner Conjecture for Almost All Points}
\author{C. Harold Horvat, Matthew Stoffregen}

\begin{document}
\maketitle

\section{The Lonely Runner Conjecture}
Consider k + 1 runners on a circular track of circumference 1, and the vector $$v = [v_1,v_2,...,v_{k+1}], v_i\in \mathbb{R}^+$$ Let d(a,b) = $||a-b||_{1}$
\\

$\textbf{Conjecture}$: $\exists  t_i \in \mathbb{R}^+ \textrm{such that }\min(d(x_i,x_j))_{i \neq j}\geq \frac{1}{k+1}.$

\section{A Solution for Almost All Points}
Suppose the components of $v$ are linearly independent. Consider the family of discrete mappings $T_i: \mathbb{S}^1 \rightarrow \mathbb{S}^1.$ defined by 
$$T_i(x)= ||x + \alpha*v_i||$$

Where

$$ \alpha \in (0,\frac{1}{\max{v_i}}) \backslash \mathbb{Q}, \frac{v_i}{\alpha} \notin \mathbb{N} \forall i.$$

 Individiually, each $T_i$ represents the irrational rotation map on $\mathbb{S}^1$ with angle $\beta = \alpha*v_i$. The trajectory of $x_i$ is dense over $\mathbb{S}^1$. Further, the product of N such rotations is dense on the torus $\mathbb{T}^k$.

	The density over $\mathbb{T}^k$ forces that every possible configuration $(x_1,x_2,...,x_k)$ is approximated arbitrarily well in finite time. Therefore, $$\forall i, \forall \epsilon, \exists n_i \in \mathbb{N}: T_{j,j\neq i}^{n_i}(0) \in (-\epsilon, \epsilon), T_i^{n_i}(0) \in (.5 - \epsilon, .5 + \epsilon)$$ We have found a time, t = $n_i * \alpha$ where runner i is lonely. Any similar bound will be surpassed, as well.

	The set of excluded points has Lebesgue measure zero. Consider that any rationally dependent velocity vector must be orthogonal to a different rational vector. Therefore it lies in a hyperplane of dimension n-1. This hyperplane has measure zero in $\mathbb{R}^n$. As each such hyperplane is identified by a rational vector orthogonal to it, there are countably many. Thus the set of all such hyperplanes, which contains the set of rationally dependent vectors, is measure zero as well. 

\section{Coprime points}
Consider the n+1-runner Lonely Runner Conjecture, set one of the runners velocities to 0.  Now observe that we need only consider the case of :
\begin{displaymath} (1/D_1,...1/D_n) \end{displaymath}
where the \( D_i \) are pairwise relatively prime integers.  (WLOG, \(D_i<D_{i+1} \)) Call such sets "best".  Then we want a time \( T \) with \( T/D_i = X_i \pmod{1}\), where the \(X_i\) are rational numbers between \( 1/(n+1), n/(n+1) \).  If we choose the \( X_i \) so that \( D_iX_i \) is an integer for all i, then we have that this is equivalent to:
\( T=X_iD_i \pmod{D_i} \) for all i.  By the Chinese Remainder Theorem, we have that there exists such a \(T<\displaystyle\prod_{i=1}^n D_i.\)  (and moreover, T is an integer).

Now, if we have maps determined by the sets of initial velocities, and some time T, then the change of position of the nth runner due to the change in initial velocities is \( \delta v * T \).  So if we choose \(\delta v\) less than 

\begin{displaymath} (1/2-1/(n+1)-1/2D_i)/\displaystyle\prod_{i=1}^n D_i \end{displaymath}

 then the set of new initial velocities is also good.  This requires some explanation.  We choose the \( X_i \) to be the closest approximation to 1/2 that the given \( D_i \) will permit.  This is reasonable, because it keeps us the farthest from going out-of-bounds on the circle.  The maximum error in approximating 1/2 by \( a/D_i) \) for some integer a, is \( 1/2D_i \).  Thus, with the given \( \delta v \), the maximum error incurred at the time T, when the original map was best, for the velopcity \( D_i \) is:
 
 \begin{displaymath}   ((1/2-1/(n+1)-1/2D_i)/\displaystyle\prod_{i=1}^n D_i)T< 1/2-1/(n+1)-1/2D_i
 \end{displaymath}
 The distance from \(X_iD_i \) to 1/n or n/(n+1), on the other hand, is at least:
 
 \begin{displaymath} 1/2-1/(n+1)-1/2D_i \end{displaymath}
 
 by our above argument, so indeed, if the \( \delta v \) is as small as demanded above, then the new set of velocities is "good".

 However, we are only concerned with the ratios between velocities.  We will find the maximum and minimum velocity ratios of the good sets generated above from the best sets.  Any set that has all of its velocity ratios lying between the maximum and minimum for each of the i in 1,...,n-1, must be good.  
 
 A velocity vector \( (x_1,...,x_n) \) is completely determined by the ratios \begin{displaymath} (x_2/x_1, x_3/x_2,...x_n/x_{n-1}) \end{displaymath}  We will now specify the \emph{most} that each of these ratios can be for a good x-vector generated from a given best vector.   
 
 Consider the ratio \( x_{i+1}/x_i \):  
 The maximum velocity for the ith component \( 1/D_i \) is
 \begin{displaymath} 1/(D_i-\delta v_i) \end{displaymath}
 
 And the minimum for the \( D_{i+1} \) component is:
 \begin{displaymath}  1/(D_{i+1}+\delta v_{i+1}) \end{displaymath}
 
 Then, if the ith ratio for some given x-vector is greater than 
 \begin{displaymath} \frac{D_{i}-\delta v_{i}}{D_{i+1}+\delta v_{i+1}} \end{displaymath}
 we are halfway there.  It is easy to see by a similar calculation that if each ratio of the x-vector is less than 
 \begin{displaymath} \frac{D_{i}+\delta v_{i}}{D_{i+1}-\delta v_{i+1}} \end{displaymath}
 then it indeed lies in the appropriate ratio zone.

 From this we can now reformulate the conjecture:
 A vector is good (with ratios \( (R_1,...R_{n-1})\), all nonzero rational)  if there exist pairwise relatively prime \( D_i \) i=1,...n-1, with:
 
 \begin{displaymath}  \frac{D_{i}-\delta v_{i}}{D_{i+1}+\delta v_{i+1}} < R_i < \frac{D_{i}+\delta v_{i}}{D_{i+1}-\delta v_{i+1}} 
 \end{displaymath}
 
 where \begin{displaymath} \delta v_i = ((1/2-1/(n+1)-1/2D_i)/\displaystyle\prod_{i=1}^n D_i \end{displaymath}
 
 \section{Small Strategy for this Style}
 Define the Quality \( Q(D_i,D_{i+1}) \) of an approximation \( D_i/D_{i+1} \) to \(R_i\) as the maximum that \begin{displaymath} \displaystyle\prod_{j=1}^n (D_j)/D_iD_{i+1} \end{displaymath} may be and to still have: 
  \begin{displaymath} \frac{D_{i}-\delta v_{i}}{D_{i+1}+\delta v_{i+1}} < R_i < \frac{D_{i}+\delta v_{i}}{D_{i+1}-\delta v_{i+1}} \end{displaymath}
  
  It is a quick observation to see that if \( \displaystyle\prod_{i=1}^n (D_j)/D_iD_{i+1}  <Q(D_i,D_{i+1}) \), for all i, then the ratio vector is good.  I think this is probably the ideal way to go about this route, if to go about it at all.  
  
  Further, it can be shown that all the Q's are almost completely determined once you have selected merely one of the D to do the approximating.  Haven't gotten to that yet.

\section{Problems}
We will generate the volume of all such points formed by the quality approach.    
For a given vector \begin{math} (P_1/D_1,...,P_n/D_n) \end{math} the (projective) volume of the set of good velocities generated by this vector is the n 
Recall:

 \begin{displaymath}  \frac{D_{i}-\delta v_{i}}{D_{i+1}+\delta v_{i+1}} < R_i < \frac{D_{i}+\delta v_{i}}{D_{i+1}-\delta v_{i+1}} 
 \end{displaymath}
 
 For all the good vectors (written as only ratios \( R_i \)).  The volume is then merely the volume of the intervals of good \( R_i\).  These are:
 
 \begin{displaymath} \displaystyle\prod_{i=1}^{n-1}(\frac{D_{i}+\delta v_{i}}{D_{i+1}-\delta v_{i+1}}-\frac{D_{i}-\delta v_{i}}{D_{i+1}+\delta v_{i+1}})=\frac{\displaystyle\prod_{i=1}^{n-1}2(\delta v_{i+1}D_i-\delta v_iD_{i+1})}{\displaystyle\prod_{i=1}^{n-1}(D_{i+1}-\delta v_{i+1})(D_{i+1}+\delta v_{i+1})}
 \end{displaymath}
 
 This then needs to be summed over all possible choices for \( P_i, D_i \).  We here order the fractions with the least \(D_i\) first, allowing negatives in the \(P_i\) only.  We will index by the number of negatives in the velocity vector.  Now,  we have the following sum:
 
\begin{displaymath}
2^{n-1}\sum_{D_1,D_2,...D_{n} coprime}\sum_{P_1,...,P_n}\frac{\displaystyle\prod_{i=1}^{n-1}(\delta v_{i+1}D_i-\delta v_iD_{i+1})}{\displaystyle\prod_{i=1}^{n-1}(D_{i+1}-\delta v_{i+1})(D_{i+1}+\delta v_{i+1})}
\end{displaymath}

\section{Going About it Backwards}
Here we will display an alternative approach to the Conjecture and provide a strategy for further efforts.  Consider once again \(R^n\), where one runner's velocity has been set to zero. 

Some definitions are in order.  Let \( ||x|| \) for \(x\in R\) be the least distance from \(||x||\) to an integer.    
An "exact case" will be any velocity vector \((v_1,...v_n)\)in \(R^n\) that has
\begin{displaymath}
sup_{t\in(0,\inf)}inf_{i=1,...,n}(||v_i*t||)=1/(n+1)
\end{displaymath}

A "pseudo-exact case" will be any velocity vector, which is not an exact case, \(v_1,...v_n\) in \(R^n\) that has    
\begin{displaymath}
sup_{t\in(0,T)}inf_{i=1,...,n}(||v_i*t||)=1/(n+1) 
\end{displaymath}
and
\begin{displaymath}
sup_{t\in(0,T+\epsilon)}inf_{i=1,...,n}(||v_i*t||)=1/(n+1)
\end{displaymath}
for some T, and some \(\epsilon>0\).  
First, observe that no rationally independent vector, say \(V=(v_1,...,v_n)\), can be pseudo-exact.  For then, we would need both \(||v_i*T||\) and \(||v_j*T||\) to be rational, for some i and j.  Then \(v_i\) and \(v_j\) would be rational multiples of each other, a contradiction.  

 Suppose there exists some rational point in the velocity space\( P=(p_1/q_1,...,p_n/q_n) \) for which the Conjecture does not hold.  Then, for any rationally independent point \( A=(a_1...,a_n) )\), let \( l \) denote the line from \(P\) to \(A\).  
It is known that for \(A\) there exists some time \(T_1\) for which \(T(A)\) has each runner at least a distance \(1/(n+1) \) from the start of the track, where \(T:R^n->R^n (v_1,v_2,...,v_n)->(||v_1*T||,...||v_N*T||)\) takes velocities to the position each of their runners has at time T.  Define the map \(F_1\) as follows:
\begin{displaymath} 
F_1=inf_{t\in(0,T_1)} sup_{i=1,...,n}(||v_i*t||)   
\end{displaymath}      
Then \(F_1(A)>1/(n+1)\).  We also verify that \(F_1\) is continuous.  Indeed, the sup and inf, taken on finite intervals, are continuous, and we merely have a composition of continuous functions.  By hypothesis, 
\begin{displaymath} 
F_1(P)<1/(n+1)
\end{displaymath}
Thus, there exists some point \(B\) on \( l \) with \(F_1(B)=1/(n+1)\).  Either it is exact, or it is pseudo-exact.  However, this point cannot be rational (by construction   
\end{document}